\documentclass[11pt,a4paper]{article}

\usepackage[margin=3cm]{geometry}
\usepackage{bbm}
\usepackage{amsmath}
\usepackage{amssymb}
\usepackage{amsfonts}
\usepackage{amsthm}
\usepackage{wasysym}
\usepackage{dsfont}
\usepackage{todonotes}
\usepackage{enumitem}
\usepackage{bm}
\usepackage{soul}
\usepackage{comment}

\usepackage{graphicx}
\usepackage{wrapfig}
\usepackage{tikz}
\usetikzlibrary{shapes.geometric}
\usepackage{orcidlink}
\usepackage{hyperref}
\usepackage{cleveref}
\usepackage[style=trad-alpha,backend=biber, backref=false, isbn=false, eprint=false]{biblatex}
\title{Computer-aided solution to the $k$-bonacci pick-up sticks problem}
\author{Julian Kern\thanks{Freie Universität Berlin, Germany.}\ \orcidlink{0000-0002-8231-0736}}
\date{}

\theoremstyle{plain}
\newtheorem*{theorem}{Theorem}

\newtheorem*{lemma}{Lemma}

\theoremstyle{definition}

\begin{document}

\maketitle

\begin{abstract}
A full solution to the recently proposed problem of determining the probability that no $k$-gon can be built from $n$ independently and uniformly chosen sticks in $[0,1]$ is proposed.
This extends the known results for triangles and quadrilaterals to general $k$-gons and offers a clearer interpretation of the connection to products of $k$-bonacci numbers.
\end{abstract}

\section{Introduction}

I recently became aware of a very nice student project that culminated in the work \cite{sudbury2025pickupsticksfibonaccifactorial} in which the following natural question, also known as \emph{pick-up sticks problem}, was investigated: Given $n$ independently drawn sticks of uniform length in $[0,1]$, what is the probability that no triangle or, more generally, no $k$-gon can be built from them.
The surprising answer for triangles is that the probability is given by the inverse Fibonacci factorial $\frac{1}{\prod_{i=1}^n F_i}$, where $F_k$ is the $k$-th Fibonacci number.
They further extend their method to quadrilaterals for which the probability involves so-called tribonacci numbers, i.e.~sequences with recursion $T_{i+3} = T_i + T_{i+1} + T_{i+2}$, and speculate that a general closed formula might exist for $k$-gons with $3\leq k\leq n$.
Note that the probability that no $n$-gon can be formed is known to be $\frac{1}{(n-1)!}$, see \cite[Corollary 5]{brokenbricks}.
Although there clearly seems to be a strong connection between the pick-up sticks problem and $k$-bonacci numbers, the proof method in \cite{sudbury2025pickupsticksfibonaccifactorial} does not give much intuition about their appearance, a problem that they leave as an open problem.

By pure chance, my discovery of that work coincided with the release of the fifth version of the LLM \emph{ChatGPT}.
So I decided to test its capabilities by challenging it with the pick-up sticks problem in the triangle case.
To my surprise, it did not only solve the problem immediately, but it did so with a technique that makes the appearance of Fibonacci numbers less mysterious.
Interestingly, the LLM was not able to solve the general $k$-gon case.
It is, however, possible to extend the technique to find the following closed formula in terms of $k$-bonacci sequences.

\begin{theorem}
Let $n\geq 3$ and $3\leq k\leq n$.
Denote for $1\leq \ell\leq k-1$ by $T_{k,\ell}$ the $(k-1)$-bonacci numbers with initial terms
\[
T_{k,\ell}(1) = \cdots = T_{k,\ell}(\ell-1) = 0\quad\text{ and }\quad T_{k,\ell}({\ell}) =\cdots = T_{k,\ell}(k-1) = 1.
\]
Then, the probability that no $k$-gon can be constructed from $n$ random lengths drawn independently and uniformly in $[0,1]$ is given by
\[
p_{n,k} = \dfrac{1}{\prod_{\ell=1}^{k-1} T_{k,\ell}(n)\cdot \prod_{i=1}^{n-k+1} T_{k,k-1}(n-i)}.
\]
Writing $S_k$ for the $(k-1)$-bonacci sequence with initial terms
\[
S_k(1) = 1\qquad\text{ and }\qquad S_k(i) = 2^{i-2}\text{ for }2\leq i\leq k-1,
\]
we can rewrite the above as
\[
p_{n,k} = \dfrac{1}{\prod_{\ell = 1}^{k-3} T_{k,\ell}(n)\cdot \prod_{i=1}^{n-k+3} S_k(i).}
\]
\end{theorem}

When $n$ is large compared to $k$, the formula is almost the inverse of a $(k-1)$-bonacci factorial.
The proof shows how the correction term emerges from the fact that the $k-1$ smallest sticks do not have to satisfy any constraint.

In the case $k = 3$, the sequence $S_k$ produces the Fibonacci numbers; in the case $k = 4$, the sequence $S_k$ specialises to produce the tribonacci sequence used in \cite{sudbury2025pickupsticksfibonaccifactorial} and we have the particularity that $T_{k,1}(n) = S_k(n) - S_k({n-2})$, simplifying the result.
In general, $T_{k,\ell}(i)$ can be expressed as a linear combination of the values of $S_k$, but the dependencies are not as simple as in the cases $k\in\{3,4\}$, which might explain the difficulty of extending the method from \cite{sudbury2025pickupsticksfibonaccifactorial} to higher values of $k$.
In the special case $k = n$, we have $T_{k,\ell}(n) = \sum_{j = \ell}^{n-1} 1 = n-\ell$, so that we recover
\[
p_{n,n} = \dfrac{1}{\left(\prod_{\ell = 1}^{n-1} T_{k,\ell}(n)\right)\cdot T_{n,n-1}(n-1)} = \dfrac{1}{(n-1)!}
\]
as expected.
For the case $k = n-1$, the computations already become more extensive: we have $T_{n-1,\ell}(n) = 2(n-1-\ell) - \delta_{\ell = 1}$ and therefore
\[
p_{n,n-1} = \dfrac{1}{(2n - 5)2^{n-3}(n-3)!}.
\]

\section{Proof of the Theorem}

As already noticed in \cite{sudbury2025pickupsticksfibonaccifactorial}, we can reformulate the problem as follows: Let $U_{(1)} < U_{(2)} < \cdots < U_{(n)}$ denote the order statistics of $n$ independent random variables with common distribution $\mathrm{Unif}([0,1])$.
Then, a $k$-gon can be constructed from these sticks if and only if
\[
U_{(i)} + \cdots + U_{(i+k-2)} \geq U_{(i+k-1)}
\]
for some $1\leq i\leq n-k+1$.
Instead of using a probabilistic argument to analyse these constraints, we take a more geometric approach.
More precisely, we use the fact that computing $p_{n,k}$ is the same as computing $n!$ times the volume of the set
\[
E_{n,k} := \left\{ u\in \mathbb R^n\;:\; \begin{cases}
    0&< u_1\\
    u_{i-1} &< u_{i} \qquad \text{ for all $2\leq i\leq k-1$}\\
    \sum_{j=i-k+1}^{i-1} u_j &< u_{i} \qquad\text{ for all $k\leq i\leq n$}\\
    u_n &\leq 1
\end{cases}\right\}.
\]
The main insight is that we can transform this complicated set into a simplex via a volume-preserving linear transformation.
More precisely, note that by setting $v_1 := u_1$, $v_i := u_i - u_{i-1}$ for $2\leq i\leq k-1$ and
\[
v_i := u_i - \sum_{j=i-k+1}^{i-1} u_j\qquad\text{ for }k\leq i\leq n,
\]
the above constraints transform into
\begin{equation}\label{eq:constraints}
\begin{cases}
    v_i > 0 &\text{ for all }1\leq i\leq n\\
    u_n \leq 1
\end{cases}\quad.
\end{equation}
It is clear that we will recover some kind of $(k-1)$-bonacci numbers when inverting the above linear transformation.
The exact relation is given in the next lemma, the proof of which we will postpone to the end of the section.

\begin{lemma}
For $1\leq m\leq k-1$, we have
\[
u_m = \sum_{j=1}^m v_j = \sum_{\ell=1}^{m} T_{k,\ell}(m)v_\ell.
\]
For $k\leq m\leq n$, the relation becomes
\[
u_m = \sum_{\ell =1}^{k-1} T_{k,\ell}(m)v_\ell + \sum_{i=1}^{m-k+1} T_{k,k-1}(m-i)v_{k-1+i}.
\]
\end{lemma}

For $m = n$, we obtain the representation
\[
u_n = \sum_{\ell=1}^{k-1} T_{k,\ell}(n)v_\ell + \sum_{i=1}^{n-k+1} T_{k,k-1}(n-i) v_{k-1+i}.
\]
As such, the constraints \eqref{eq:constraints} define a simplex with vertices 
\[
0, \frac{1}{T_{k,1}(n)}\mathfrak{e}_1, \dots, \frac{1}{T_{k,k-1}(n)}\mathfrak{e}_{k-1}, \frac{1}{T_{k,k-1}(n-1)}\mathfrak{e}_k, \dots, \frac{1}{T_{k,k-1}(k-1)}\mathfrak{e}_n,
\]
where $(\mathfrak{e}_1,\dots, \mathfrak{e}_n)$ denotes the standard basis of $\mathbb R^n$.
With the vertices forming the columns of a diagonal matrix, we deduce that its volume is given by
\[
\frac{1}{n!}\cdot \left(\prod_{\ell = 1}^{k-1} \dfrac{1}{T_{k,\ell}(n)}\right)\cdot \left(\prod_{i=1}^{n-k+1} \dfrac{1}{T_{k,k-1}(n-i)}\right),
\]
which concludes the proof.

\begin{proof}[Proof of the Lemma]
For $1\leq m\leq k-1$, this follows from the fact that
\[
u_m = u_1 + \sum_{j=2}^m (u_j - u_{j-1}) = \sum_{j=1}^m v_j.
\]
For $m\geq k$, we first note that if we have the representation $u_p = \sum_{j=1}^p a_{p,j}v_p$ for $p\leq m-1$, then
\[
u_m = v_m + \sum_{p=1}^{k-1}u_{m-p} = \sum_{p=1}^{k-1}\sum_{j=1}^{m-p} a_{m-p,j}v_j = \sum_{j=1}^{m-1} \left(\sum_{p = 1}^{k-1} a_{m-p,j}\right) v_j,
\]
where we set $a_{m-p,j} = 0$ whenever $m-p < j$.
This means in particular that $a_{m,m} = 1$ and that the coefficients satisfy the recurrence relation of $(k-1)$-bonacci numbers.
To identify the initial values, we note that for $j\leq m\leq k-1$, we have by the first point $a_{m,j} = 1$.
For $j \geq m+1$, we will always have $a_{m,j} = 0$.
In particular, for $\ell\leq k-1$, we will have $a_{1,\ell} = \cdots = a_{\ell-1, \ell} = 0$ and $a_{\ell,\ell} = \cdots = a_{k-1, \ell} = 1$ which uniquely characterises the sequence $T_{k,\ell}$.
For $k\leq j\leq n$, we automatically have $a_{m,j} = 0$ for all $1\leq m\leq j-1$.
Furthermore, recall that $a_{m,m} = 1$ by the above computation.
But this uniquely identifies $m\mapsto a_{m + j -k+1, j}$ to be the sequences $T_{k,k-1}$.
\end{proof}

\printbibliography

@misc{sudbury2025pickupsticksfibonaccifactorial,
      title={Pick-up Sticks and the {F}ibonacci Factorial}, 
      author={A. Sudbury and A. Sun and D. Treeby and E. Wang},
      year={2025},
      eprint={2504.19911},
      archivePrefix={arXiv},
      primaryClass={math.PR},
      url={https://arxiv.org/abs/2504.19911}, 
}

@article{brokenbricks,
 ISSN = {0025570X, 19300980},
 URL = {https://www.jstor.org/stable/48665725},
 author = {T. Petersen and B. Tenner},
 journal = {Mathematics Magazine},
 number = {3},
 pages = {pp. 175--185},
 publisher = {[Mathematical Association of America, Taylor & Francis, Ltd.]},
 title = {Broken Bricks and the Pick-up Sticks Problem},
 volume = {93},
 year = {2020}
}

\end{document}